\documentclass{elsart}
\usepackage{amssymb}
\usepackage{graphicx}

\begin{document}

\begin{frontmatter}

\title{Towards~Theory~of~Piecewise~Linear 
Dynamical~Systems\thanksref{label1}}
\thanks[label1]{This work is supported by the Netherlands Organization for Scientific
Research (NWO).}

\author{Valery A. Gaiko\thanksref{label2}}
\ead{vlrgk@yahoo.com}
\thanks[label2]{The first author is grateful to the Delft Institute of Applied Mathematics 
of TU~Delft for hospitality during his stay at the University in the period of 
October, 2007~-- March, 2008.}

\address{Belarusian State University of Informatics and Radioelectronics,
Department~of~Mathematics,~L.\,Beda~Str.\,6--4,~Minsk~220040,~Belarus}

\author{Wim T. van Horssen}
\address{Delft University of Technology, Delft Institute of Applied Mathematics,
Mekelweg~4,~2628~CD~Delft,~The~Netherlands}

\begin{abstract}
In this paper, we consider a planar dynamical system with a piecewise linear function containing an arbitrary number (but finite) of dropping sections and approximating some continuous nonlinear function. Studying all possible local and global bifurcations of its limit cycles, we prove that such a piecewise linear dynamical system with $k$ dropping sections and $2k+1$ singular points can have at most $k+2$ limit cycles, $k+1$ of which surround the foci one by one and the last, $(k+2)$-th, limit cycle surrounds all of the singular points of this system.
\end{abstract}

\end{frontmatter}

\section{Introduction}

The paper is based on the applications of Bifurcation Theory developed by Andronov, Arnold, Thom, Whitney, Zeeman et al. and 
can be used for modeling problems, where system parameters play a certain role in various bifurcations. The theoretical studies of bifurcations deal with so-called universal problems. This means that sufficiently many parameters are available for universality of generic families of dynamical systems in the context at hand, under a relevant equivalence relation. This has led to the classification of generic, local bifurcations. In many applications, models have a given number of parameters. Moreover, the bifurcation analysis, taking place in the product of phase space and parameter space, is not restricted to local features only. On the contrary, often the interest is the global organization of the parameter space regarding bifurcations which can be both local and global. 

This paper deals with so-called sewed dynamical systems, i.e., with systems for which the domain of definition is divided into sub-domains where different analytical systems are defined. The trajectories of these partial systems are sewed in one way or another on the boundaries of the sub-domains. Such systems have some typical features, namely: 1)~the system sewing is immediate from the physical meaning of the problem under consideration; 2)~the system is piecewise linear, i.e., the partial systems from which it is sewed are linear systems; 3)~on the line of sewing, a point map (a first return function) is defined, what allows to determine the character of the system under consideration. 

Piecewise linear dynamical systems always contain some parameters and, under the variation of these parameters, the qualitative behavior of the systems can obviously change. We will consider the simplest bifurcations possible in the sewed systems when the sewing lines are unchanged under the parameter variations. It is natural to consider the following bifurcations which are similar to the simplest bifurcations of continuous dynamical systems: 1)~the bifurcation of a singular point of focus type; 2)~the bifurcation of an immovable point of focus type, a quasi-focus; 3)~the bifurcation of a sewed limit cycle; 4)~the bifurcation of a sewed separatrix going from a saddle to another saddle (the saddles can be both sewed and unsewed); 5)~the bifurcation of a separatrix of a saddle-shaped singular point (sewed or unsewed) going from a saddle-shaped singular point to another such point or to a saddle (sewed or unsewed); 6)~the bifurcation of a sewed saddle-node; 7)~the bifurcation of a sewed separatrix of a saddle-node (sewed or unsewed) going out of the saddle-node and going back to it. Besides, some specific bifurcations can occur in sewed systems. Since in such systems, for example, arches of attraction or repulsion composed of immovable points can be similar to singular points, some bifurcations which are similar to the generation of a limit cycle from a focus can occur in the corresponding constructions.

Piecewise linear systems have many applications in science and engineering. Special cases of such systems provide mathematical models for mechanical systems with Coulomb friction, for valve oscillators with a discontinuous characteristic, for direct control systems with a two-point relay mechanism, for planar dynamical systems modeling neural activity, etc. Despite their simple structure and relevance to the applications, there is, to the best of our knowledge, no complete study of their dynamical properties. In most existing papers, which deal with planar dynamical systems with piecewise linear right-hand sides, either the systems considered are continuous or only particular cases are investigated. The first analytical results on such systems go back to Andronov, Vitt, and Khaikin in the 1930s (see \cite{1}). The  existence and non-existence of an asymptotically stable periodic solution (limit cycle) of a piecewise linear system can be comparatively easily proved. However, for example, periodic solutions with sliding motion are of great importance to the applications. In~particular, they describe the so-called stick-slip oscillations which appear in mechanical systems with dry friction. 

The main objective of the present paper is to provide a complete analysis of the dynamical properties of piecewise linear systems, their dependence on the system parameters studying, first of all, their limit cycle bifurcations. There are several ways to investigate the qualitative dynamics of such systems \cite{7}. There are also numerous methods and good results on studying limit cycles. However, the most important impulse to their studying was given by the introduction of ideas coming from Bifurcation Theory \cite{2}, \cite{3}, \cite{8}, \cite{13}. We know three principal limit cycle bifurcations which we will study~\cite{8}: 1)~the Andronov-Hopf bifurcation (from a singular point of centre or focus type); 2)~the separatrix cycle bifurcation (from a~homoclinic or heteroclinic orbit); 3)~the multiple limit cycle bifurcation. All of these local bifurcations will be globally connected and will be applied to the qualitative analysis of piecewise linear dynamical systems.

\section{Preliminaries}

In this paper, geometric aspects of Bifurcation Theory are used and  developed. It gives a global approach to the qualitative analysis and helps to combine all other approaches, their methods and results. First of all, the two-isocline method which was developed by Erugin is used~\cite{8}. An isocline portrait is the most natural construction for a polynomial equation. It is sufficient to have only two isoclines (of zero and infinity) to obtain principal information on the original polynomial system, because these two isoclines are right-hand sides of the system. Geometric properties of isoclines (conics, cubics, quadrics, etc.) are well-known, and all isoclines portraits can be easily constructed. By means of them, all topologically different qualitative pictures of integral curves to within a number of limit cycles and distinguishing center and focus can be obtained. Thus, it is possible to carry out a rough topological classification of the phase portraits for the polynomial dynamical systems and for the corresponding piecewise linear systems. It is the first application of Erugin's method. After studying contact and rotation properties of isoclines, the simplest (canonical) systems containing limit cycles can be also constructed. Two groups of the parameters can be distinguished in such systems: static and dynamic. Static parameters determine the behavior of phase trajectories in principle, since they control the number, position and character of singular points in a finite part of the plane (finite singularities). The parameters from the first group determine also a possible behavior of separatrices and singular points at infinity (infinite singularities) under the variation of the parameters from the second group. The dynamic parameters are rotation parameters. They do not change the number, position and index of the finite singularities, but only involve the vector field in the directional rotation. The rotation parameters allow to control the infinite singularities, the behavior of limit cycles and separatrices. The cyclicity of singular points and separatrix cycles, the behavior of semi-stable and other multiple limit cycles are controlled by these parameters as well. Therefore, by means of the rotation parameters, it is possible to control all limit cycle bifurcations and to solve the most complicated problems of the qualitative theory of dynamical systems (both continuous and piecewise linear).

In \cite{8} some complete results on continuous quadratic systems have been presented, and some preliminary results on generalizing geometric ideas and bifurcation methods to cubic dynamical systems have been obtained. So, in \cite{11},  a canonical cubic dynamical system of Kukles-type was constructed and the global qualitative analysis of its special case corresponding to a generalized Li\'{e}nard equation was given. In particular, it was proved that the foci of such a Li\'{e}nard system could be at most of second order and that the system could have at least three limit cycles in the whole phase plane. Moreover, unlike all previous works on the Kukles-type systems, by means of arbitrary (including as large as possible) field rotation parameters of the canonical system, the global bifurcations of its limit and separatrix cycles were studied. As a result, the classification of all possible types of separatrix cycles was obtained and all possible distributions of limit cycles were found for the generalized Li\'{e}nard system. In \cite{5}, the global qualitative analysis of centrally symmetric cubic systems which are used as learning models of planar neural networks was established. All of these results can be generalized to the corresponding piecewise linear Li\'{e}nard-type dynamical systems.

Such systems have been considered in several paper. For example, Bautin \cite{3} studied a dynamical system which was used in radio-engineering for describing tunnel diode circuits, where the nonlinear function was approximated by a piecewise linear function composed of three linear pieces. Giannakopoulos and Pliete \cite{12} generalized that results for a piecewise linear system with a line of discontinuity modeling stick-slip oscillations which appeared in mechanical systems with dry friction.  In this paper, the obtained results are generalized to an arbitrary piecewise linear Li\'{e}nard-type system, where the polynomial is approximated by a piecewise linear function composed of an arbitrary (but finite) number of linear pieces. In particular, the following bifurcations are studied for such a system: 1)~the bifurcation of singular points (division of the parameter space according to the number and character of singular points, stability of singular points lying on the lines of sewing, generation of limit cycles from singular points of focus type under transferring the points through the lines of sewing, generation of limit cycles (hyperbolic and semi-stable) from the boundaries of the domains filled by closed trajectories); 2)~the separatrix bifurcation (location of the bifurcation curves for the separatrix loops, stability of the separatrix loops); 3)~the bifurcation of multiple limit cycles (location of the bifurcation curves for the limit cycles of various multiplicity, qualitative structure of the division in the phase plane).

\section{Main Results}

Consider the system
\begin{equation}
	\label{E1}
\dot{x}=y-\varphi(x), \
\dot{y}=\beta-\alpha x-y, \
	\alpha>0, \ \beta>0,
\end{equation}
where $\varphi(x)$ is a piecewise linear function containing $k$ dropping sections and approximating some continuous nonlinear function. The line $\beta-\alpha x-y=0$ and the curve $y=\varphi(x)$ can be considered as the isoclines of zero and infinity, respectively, for the corresponding equation. Such systems and equations may occur, for example, when tunnel diode circuits and some other problems are studied (see \cite{1}, \cite{3}, \cite{4}, \cite{7}, \cite{13}).

Suppose that the ascending sections of (\ref{E1}) have an inclination $k_{1}>0$ and the descending (dropping) sections have an inclination $k_{2}<0.$ Then the phase plane of (\ref{E1}) can be divided onto $2k+1$ parts in every of which (\ref{E1}) is a linear system: the ascending sections are in $k+1$ strip regions $(I, III ,V,\ldots, 2K+1)$ and the descending sections are in other $k$ such regions $(II, IV, VI,\ldots, 2K).$ The parameters $k_{1},$ $k_{2},$ and also $\alpha$ can be considered as rotation parameters for the sewed vector field of (\ref{E1}) (see \cite{4}, \cite{8}).

System (\ref{E1}) can have an odd number of rough singular points: $1,  3,  5,\ldots, 2k+1.$ If (\ref{E1}) has the only singular point, this point will be always an antisaddle (center, focus or node). A focus (node) will be always stable in odd regions and unstable in even regions if $k_{2}>1.$ If system (\ref{E1}) has $2k+1$ singularities, then $k$ of them are saddles (they are in even regions) and $k+1$ others are antisaddles (foci or nodes) which are always stable (they are in odd regions).  The pieces of the straight lines  $\beta=x_{2i-1}\alpha+y_{2i-1}$ and $\beta=x_{2i}\alpha+y_{2i}$\linebreak $(i=1,2,\ldots,k),$ where $(x_{2i-1},y_{2i-1})$ and $(x_{2i},y_{2i})$ are the coordinates of the upper and lower corner points of the curve $\varphi(x),$ respectively, form a discriminant curve separating the domains in the plane $(\alpha,\beta),$ where $\alpha\leq k_{2},$ with different numbers of singular points. The points of the discriminant curve correspond to the sewed singularities of saddle-focus or saddle-node type $(\alpha<k_{2})$ and its corner points correspond to the unstable equilibrium segments $(\alpha=k_{2})$ which coincide with the dropping sections of the curve $y=\varphi(x).$ In the case when $k_{2}<1,$ closed trajectories cannot exist and only bifurcations of singular points are possible in system (\ref{E1}). Therefore, we will consider further only the case when $k_{2}>1$ and $(k_{1}-1)^{2}<4k_{2}$ giving various bifurcations and, first of all, the bifurcations of limit cycles. 

In Section 4, studying local bifurcations of limit cycles, we prove that: 1)~the only limit cycle can appear in system (\ref{E1}) under a displacement of a singular point of focus type through the sewing line of (\ref{E1}), and this limit cycle will be unique in the corresponding pair of the strip regions; 2)~at most two limit cycles can be generated by the boundary of the domain filled by closed trajectories of system (\ref{E1}), and these limit cycles can be only outside the boundary; 3)~the equilibrium segment (the dropping section) of (\ref{E1}) can generate at most two unstable limit cycles surrounging two stable foci one by one which appear simultaneously from the ends of the segment; 4)~a nondegenerate separatrix cycle of system (\ref{E1}) with $2k+1$ simple singular points can generate at most $k+1$ unstable limit cycles inside its loops (digons) or the only unstable limit cycle outside it.
 In Section 5, using the obtained results and applying the Wintner--Perko termination principle for sewed multiple limit cycles, we sudy global bifurcations of limit cycles and prove that system (\ref{E1}) with $k$ dropping sections and $2k+1$ singular points can have at most $k+2$ limit cycles, $k+1$ of which surround the foci one by one and the last, $(k+2)$-th, limit cycle surrounds all of the singular points of (\ref{E1}).

\section{Local Bifurcations of Limit Cycles}

Following \cite{3}, we begin with studying local bifurcations of limit cycles. The limit cycle of system (\ref{E1}) will be called \emph{small} if it belongs to at most  two adjoining regions; the cycle will be called \emph{big} if it belongs to at least three adjoining regions.

\subsection{The Limit Cycle Bifurcation from a Singular Point of Focus Type under a Displacement of the Singular Point through the Line of Sewing}

\noindent \textbf{Lemma 4.1.}
\emph{Under a displacement of a singular point of focus type through the line of sewing of system (\ref{E1}), at most one limit cycle can appear in this system, and this limit cycle will be unique in the corresponding pair of the strip regions of (\ref{E1}).}
\medskip

\noindent \textbf{Proof.} Let us study first stability of a singular point on the line of sewing. Suppose that the straight line $\beta-\alpha x-y=0$ passes through the corner point $(x_{1},y_{1})$ of the curve $y=\varphi(x)$ on the boundary of regions $I,$ $II$ and that $\alpha>(k_{2}+1)^{2}/4.$ Then region $I$ $(II)$ will be filled by pieces of trajectories of a stable (unstable) focus of system (\ref{E1}).

Introduce positive coordinates $S_{0}$  (lower $(x_{1},y_{1}))$ and $S_{1}$  (upper $(x_{1},y_{1}))$ on the line of sewing of regions $I$ and $II;$ $S_{2}$  (lower $(x_{2},y_{2}))$ and $S_{3}$  (upper $(x_{2},y_{2}))$ on the line of sewing of regions $II$ and $III,$ etc. The maps $S_{0}\rightarrow S_{1}$ along the trajectories of region $I$ and $S_{1}\rightarrow S_{0}$ along the trajectories of region $II$ are written as follows:
\begin{equation}
	\label{E2}
S_{1}=S_{0}e^{ \pi\sigma_{1}/\omega_{1}}, \quad
	 \bar{S}_{0}=S_{1}e^{ \pi\sigma_{2}/\omega_{2}},
\end{equation}
where $\sigma_{i},$ $\omega_{i}$ $(i=1,2)$ are the real and imaginary parts of the roots of the characteristic equation for a singular point of regions $I,$ $II,$ respectively.

The singular point $(x_{1},y_{1})$ will be a sewed center $(\bar{S}_{0}\!=\!S_{0})$ iff $\sigma_{1}/\omega_{1}+\sigma_{2}/\omega_{2}\!=\!0,$ i.\,e., when
$\alpha=\alpha^{*}\equiv(1-k_{1}/k_{2})/(k_{2}-k_{1}+2).$ The sewed focus $(x_{1},y_{1})$ will be stable  $(\bar{S}_{0}<S_{0})$  when $\alpha>\alpha^{*}$ and unstable $(\bar{S}_{0}>S_{0})$ when $\alpha<\alpha^{*}.$

Consider the return map $S_{0}\rightarrow\bar{S}_{0}$ along the trajectories of regions $I$ and $II.$ For region $I,$ we will have
\begin{equation}
	\label{E3}
		\begin{array}{ll}
S_{0}&=\displaystyle\frac{\delta_{0}}{\sin\omega_{1}\tau_{1}}
	(\omega_{1}\cos\omega_{1}\tau_{1}
	-\sigma_{1}\sin\omega_{1}\tau_{1}
	-\omega_{1}e^{-\sigma_{1}\tau_{1}})
	 \equiv\delta_{0}\zeta(\tau_{1}),
		\\[4mm]
S_{1}&=\displaystyle\frac{\delta_{0}}{\sin\omega_{1}\tau_{1}}
	(\omega_{1}\cos\omega_{1}\tau_{1}
	+\sigma_{1}\sin\omega_{1}\tau_{1}
	-\omega_{1} e^{\sigma_{1}\tau_{1}})
	 \equiv\delta_{0}\chi(\tau_{1}),
		\end{array}
\end{equation}
where $\delta_{0}$ is the distance from the boundary of regions $I,$ $II$ to the singular point; $\zeta$ and $\chi$ are monotonic functions. The return map along the trajectories of region $II$ has a similar form.

Calculation of the first derivative for the return map gives
\begin{equation}
	\label{E4}
\frac{\textrm{d}\bar{S}_{0}}{\textrm{d}S_{0}}=\frac{S_{0}}{\bar{S}_{0}}\,
	e^{2( \sigma_{1}\tau_{1}+\sigma_{2}\tau_{2})},
\end{equation}
where $\tau_{i}$ $(i=1,2)$ is motion time along the trajectories of regions $I,$ $II,$ respectively; $\sigma_{i}=(1+k_{i})/2$ $(i=1,2).$

Suppose that the singular point of system (\ref{E1}) is in region $I.$ Then, along a periodic solution of (\ref{E1}) (when $\bar{S}_{0}=S_{0})$ and for increasing $S_{0},$ time $\tau_{1}$ decreases (to the value $\pi/\omega_{1}),$ time $\tau_{2}$ increases (to the value $\pi/\omega_{2}),$ and derivative (\ref{E4}) also increases. Therefore, at most two intersection points of the graph of the return map $S_{0}\rightarrow\bar{S}_{0}$ with the bisectrix of the plane $(S_{0},\bar{S}_{0})$ can exist; moreover, the fixed point with the smaller coordinate must be stable and the fixed point with the grater coordinate must be unstable. Since, under our assummption, the singular point is in region $I$ and it is a stable focus which cannot be surrounded by a stable limit cycle, then at most one limit cycle can exist in regions $I$ and $II;$ moreover, this limit cycle can be only unstable.

Suppose that the singular point of system (\ref{E1}) is in region $II.$ Then, for increasing $S_{0},$ time $\tau_{1}$ increases and time $\tau_{2}$ decreases. In a similar way, we can show that at most one (stable) limit cycle can exist in this case as well.

Let the straight line $\beta-\alpha x-y=0$ pass through the upper corner point $(x_{1},y_{1})$ of the curve $y=\varphi(x)$ again and suppose that $\alpha>\alpha^{*}$ (the case when $\alpha<\alpha^{*}$ is considered absolutely similarly). The sewed focus is stable in this case. The trajectory passing through the lower corner point $(x_{2},y_{2}),$ by (\ref{E2}), tends to the singular point for $t\rightarrow+\infty.$ This trajectory remains a spiral under a small displacement of the line $\beta-\alpha x-y=0.$ If, under such 
a displacement, the singular point gets into region $II,$ then it becomes unstable and, hence, at least one stable limit cycle appears in system (\ref{E1}). As follows from above, this cycle is unique. If, under a displacement, the singular point gets into region $I,$ obviously, limit cycles cannot appear in system (\ref{E1}). The same results can be obtained for regions  $III$ and $IV,\ldots,$  $2K-1$ and $2K.$ The lemma is proved.$\qquad\square$

\subsection{The Limit Cycle Bifurcation from the Boundary of the Domain Filled by Closed Trajectories}

\noindent \textbf{Lemma 4.2.}
\emph{The boundary of the domain filled by closed trajectories of system (\ref{E1}) can generate at most two limit cycles, and these limit cycles can be only outside the boundary.}
\medskip

\noindent \textbf{Proof.} Consider now the map $\bar{S}_{0}=f(S_{0})$ sewed of two pieces: $\bar{S}_{0}=\xi(S_{0})$ along the trajectories in regions $I,$ $II,\ldots,$ $2K$ and $\bar{S}_{0}=\psi(S_{0})$ along the trajectories in all regions, $I,$ $II,\ldots,$ $2K,$ $2K+1.$ The map $S_{0}\rightarrow S_{1}$ in region $I$ is given by (\ref{E3}). The maps $S_{1}\rightarrow S_{3},$ $S_{3}\rightarrow S_{5},\ldots,$  $S_{2k-1}\rightarrow S_{2k-2}$ $(S_{2k-1}\rightarrow S_{2k+1},$ $S_{2k+1}\rightarrow S_{2k},$ $S_{2k}\rightarrow S_{2k-2}),$ $S_{2k-2}\rightarrow S_{2k-4},\ldots,$ $S_{2}\rightarrow S_{0}$ have similar forms.

The derivatives for the functions $\xi(S_{0}),$ $\psi(S_{0})$ are given by the following expressions, respectively:
\begin{equation}
	\label{E5}
		\begin{array}{ll}
\displaystyle\frac{\textrm{d}\bar{S}_{0}}{\textrm{d}S_{0}}
	=\displaystyle\frac{S_{0}}{\bar{S}_{0}}\,
		e^{2(\sigma_{1}(\tau_{1}+\tau_{3}^{+}+\tau_{3}^{-}+\ldots+\tau_{2k-1})
		+\sigma_{2}(\tau_{2}^{+}+\tau_{2}^{-}+\ldots+\tau_{2k-2}^{+}+\tau_{2k-2}^{-}))},
		\end{array}
\end{equation}
\begin{equation}
	\label{E6}
		\begin{array}{ll}
\displaystyle\frac{\textrm{d}\bar{S}_{0}}{\textrm{d}S_{0}}
	=\displaystyle\frac{S_{0}}{\bar{S}_{0}}\,
		e^{2(\sigma_{1}(\tau_{1}+\tau_{3}^{+}+\tau_{3}^{-}+\ldots+\tau_{2k+1})
		+\sigma_{2}(\tau_{2}^{+}+\tau_{2}^{-}+\ldots+\tau_{2k}^{+}+\tau_{2k}^{-}))},
		\end{array}
\end{equation}
where $\tau_{1}, $ $\tau_{2k-1},$ $\tau_{2k+1}$ are motion times in regions $I,$ $2K-1,$ $2K+1$ and $\tau_{2i}^{+}$  $(\tau_{2i}^{-}),$ $\tau_{2i+1}^{+}$ $(\tau_{2i+1}^{-}),$ $i=1,2,\ldots,k,$ are  motion times in the upper (lower) parts of regions $II,III,\ldots,2K,$ respectively.

Let $S_{0}=S_{0}^{*}$ be the boundary value separating the intervals of definition of the maps $\xi(S_{0})$ and $\psi(S_{0}).$ 
By (\ref{E5}) and (\ref{E6}), the function $\bar{S}_{0}=f(S_{0})$ is differentiable at the point $S_{0}^{*}.$ 

Suppose now that the straight line $\beta-\alpha x-y=0$ passes through the corner point $(x_{1},y_{1})$ and that $\alpha=\alpha^{*}.$ Let us show that system (\ref{E1}) has no limit cycle in this case.

The return map in the plane $(S_{0},\bar{S}_{0})$ is sewed of two pieces: the bisectrix segment $\bar{S}_{0}=S_{0}<S_{0}^{*}$ and the curve $\bar{S}_{0}=\psi(S_{0}).$ Since the function\linebreak $\bar{S}_{0}=f(S_{0})$ is differentiable at the sewing point $S_{0}^{*},$ then $\textrm{d}\bar{S}_{0}/\textrm{d}S_{0}=1$ for $\alpha=\alpha^{*}$ (besides, if follows from (\ref{E6}) that
$\textrm{d}^{2}\bar{S}_{0}/\textrm{d}S_{0}^{2}<0).$ For increasing $S_{0}$ (from the value $S_{0}^{*}),$ the exponent degree in (\ref{E6}) will monotonically decrease at the sewing point (from the zero value). It is clear that the curve  $\bar{S}_{0}=\psi(S_{0})$ cannot have other points of intersection (or tangency) with the bisectrix in addition to $S_{0}=S_{0}^{*}.$ For $S_{0}>S_{0}^{*},$ the curve $\bar{S}_{0}=\psi(S_{0})$ lies below the bisectrix. Therefore, the spirals sewed of the trajectories of (\ref{E1}) in regions $I,II,\ldots,2K,2K+1$ tend to the boundary $S_{0}$ in regions $I,II,\ldots,2K.$ 

For a small variation of the parameters $\alpha$ and $\beta,$ the return map graph of the changed system will be in a small neighborhood of the return map graph of the original system. If we move along the half-line $L_{1}=0$ $(L_{1}\equiv\beta-\alpha x_{1}-y_{1,}$ $\alpha>k_{2})$ from the value $\alpha=\alpha^{*}$ to the side of decreasing $\alpha,$ then the return map will be either the straight line passing through the origin above the bisectrix (for $S_{0}<S_{0}^{*})$ or the curve $\bar{S}_{0}=\psi(S_{0})$ intersecting the bisectrix one time (for $S_{0}>S_{0}^{*}).$ In this case, the only stable limit cycle occurs from the boundary of the domain filled by closed trajectories. For further decreasing $\beta,$ the return map graph will move from the origin along the $S_{0}$-axis and will intersect the bisectrix two times (an unstable limit cycle occurs from the focus under its displacement through the line of sewing). If we move along the half-line to the side of increasing $\alpha$ from the value $\alpha=\alpha^{*},$ decreasing $\beta$ after that, then the return map graph will be completely below the bisectrix. It follows from continuity and differentiability of the return map that there exists a point $(\alpha^{*},\beta^{*})$ on the half-line $L_{1}=0$ at which a bifurcation curve of double limit cycles is created and then extended into some half-neighborhood of this point (below the half-line) bounding a domain with two limit cycles in the parameter space. 

The points of this bifurcation curve correspond to the tangency points of the return map graph with the bisectrix of the plane $(S_{0},\bar{S}_{0}).$ Such tangency is impossible for $S_{0}<S_{0}^{*}.$ Therefore, a double limit cycle can appear only from the boundary of the domain filled by closed trajectories of (\ref{E1}) when $S_{0}=S_{0}^{*}.$ Thus, we have proved that at most two limit cycles can be generated by the boundary and that these two limit cycles can be only outside the boundary.$\qquad\square$

\subsection{The Limit Cycle Bifurcation from the Ends of the Equilibrium Segment}

\noindent \textbf{Lemma 4.3.}
\emph{The equilibrium segment (the dropping section) of system (\ref{E1}) can generate at most two unstable limit cycles surrounging two stable foci one by one which appear simultaneously from the ends of the segment.}
\medskip

\noindent \textbf{Proof.} Suppose that a part of the straight line $\beta-\alpha x-y=0$ coincides with a dropping section of (\ref{E1}), for example, with the first one $(\alpha=k_{2}).$ The dropping section of (\ref{E1}) will be an unstable equilibrium segment and regions $I,$ $II$ (because of the condition $(k_{1}-1)^{2}<4k_{2})$ will be filled by trajectories of the stable foci. It is easy to obtain an explicit expression for the map of the half-line $S_{0}$ into itself:
\begin{equation}
	\label{E7}
\bar{S}_{0}=S_{0}\,
	e^{2\pi\sigma_{1}/\omega_{1}}+\delta(k_{2}-1)(1+e^{\pi\sigma_{1}/\omega_{1}}),
\end{equation}
where $\delta$ is the width of regions $II.$ This map has the only stable fixed point.

Rotate the straight line $\beta-\alpha x-y=0$ counterclockwise around some point on the dropping section. The equilibrium segment will be destroyed and three singular points will appear: a saddle in region $II$ and two stable foci in regions $I$ and $III.$ Let $\alpha=k_{2}-\varepsilon,$ where $\varepsilon$ is positive and sufficiently small. Being restricted by the first degree of $\varepsilon,$ we can obtain the following separatrix slopes: $-1+\varepsilon/(k_{2}-1)$ for the $\alpha$-separatrix and $-k_{2}-\varepsilon/(k_{2}-1)$ for the $\omega$-separatrix.

For $\alpha=k_{2},$ the trajectories going out the point, where a saddle appears when $\varepsilon\neq\,0,$ tend to a big stable limit cycle surrounding all of the singular points; for $\varepsilon>0,$ the $\alpha$-separatrices of the saddle from region $II$ lie in a small neighborhood of the trajectories going out the same point for $\varepsilon=0$ and, hence, they also tend to the same limit cycle. Therefore, the $\omega$-separatrices can just go out unstable limit cycles lying in regions $I-II,$ $II-III$ and surrounding stable foci appearing in regions $I$ and $II$ under the rotation of the straight line $\beta-\alpha x-y=0.$

Thus, under the rotation of this straight line, unstable limit cycles surrounding stable foci appear from the ends of the equilibrium segment (the cycles and the foci appear simultaneously). There is the only limit cycle in a neighborhood of each of the foci. This follows from that the derivative constructed by using the trajectories of the saddle from region $II$ will be also given by expression (\ref{E4}), with the only difference: for increasing $S_{0},$ $\tau_{2}\rightarrow+\infty.$ The same result can be obtained for all of the other dropping sections. This completes the proof of the lemma.$\qquad\square$

\subsection{The Separatrix Cycle Bifurcations}

The simplest type of separatrix cycles of (\ref{E1}) is a so-called eight-loop formed by two ordinary saddle loops. In the case of $2k+1$ simple singular points, a separatrix cycle can contain $k+1$ saddle loops, the first and the last of which are ordinary loops with one rough saddle on each and the $k-1$ others are separatrix digons with two rough saddles on each. Such a separatrix cycle will be called \emph{nondegenerate}. In the cases when the straight line $\beta-\alpha x-y=0$ passes through the corner points of the curve $y=\varphi(x),$ we will have \emph{degenerate} separatrix cycles of lips-type containing one or two sewed saddle-nodes. Let us consider the bifurcations of nondegenerate separatrix cycles. To consider the cases of degenerate separatrix cycles, we have either to generalize the results of \cite{6} or to apply the Wintner--Perko termination principle \cite{8},~\cite{14} which will be done in the next section.

\medskip
\noindent \textbf{Lemma 4.4.}
\emph{A nondegenerate separatrix cycle of system (\ref{E1}) with $2k+1$ simple singular points can generate at most $k+1$ small unstable limit cycles inside its loops (digons) or the only big unstable limit cycle outside it.}
\medskip

\noindent \textbf{Proof.} Let us show first that the bifurcations of a separatrix loop do not depend on the parameter $\beta.$ Fix the parameters $k_{1},$ $k_{2},$ $\alpha$ and suppose that the line $\beta-\alpha x-y=0$ passes through the upper corner point $(x_{1},y_{1})$ of the curve $y=\varphi(x)$ for some value $\beta=\beta_{0}.$ Varying $\beta$ by the value $\kappa$ $(\kappa=\beta_{0}-\beta),$ we will show that a separatrix loop cannot be formed under such a variation. Let $S'_{0}$ and $S'_{1}$ be the segments intercepted on the boundary of regions $I$ and $II;$ $S_{0}$ and $S_{1}$ be the coordinates for return map (\ref{E3}) on the same boundary. It follows from (\ref{E3}) that
\begin{equation}
	\label{E8}
	S_{1}=\delta_{0}\chi\zeta^{-1}(S_{0}/\delta_{0}),
\end{equation}
where $\zeta^{-1}$ is the inverse function for $\zeta.$ The values $\sigma_{1}$ and $\omega_{1}$ do not depend on $\beta;$ therefore, the functions $\chi$ and $\zeta$ also do not depend on it.

Since $y=\varphi(x)$ is a piesewise linear function, then, under the variation of $\beta,$ the values $S'_{0},$ $S'_{1},$ and $\delta_{0}$ will be proportional to $\kappa:$
\begin{equation}
	\label{E9}
	S'_{0}=\gamma_{0}\kappa,\quad S'_{1}=\gamma_{2}\kappa,\quad \delta_{0}=\gamma_{1}\kappa.
\end{equation}

Sewing the trajectories on the boundary of regions $I$ and $II$ (assuming $S'_{0}~=~S_{0}),$ we can find
from (\ref{E8}), (\ref{E9}):
\begin{equation}
	\label{E10}
	S_{1}=\gamma_{1}\kappa\chi\zeta^{-1}(\gamma_{0}/\gamma_{1})\equiv\gamma_{3}\kappa,
\end{equation}
\begin{equation}
	\label{E11}
	S_{1}/S'_{1}=\gamma_{3}/\gamma_{2}=\textrm{const}.
\end{equation}

Thus, under the fixed $k_{1},$ $k_{2},$ and $\alpha,$ the values $S_{1}$ and $S'_{1}$ have a constant ratio; hence, a separatrix loop $(S_{1}= S'_{1} )$  cannot be formed under a variation of the parameter $\beta.$ The same conclusion can be made for any other type of the separatrix cycles of (\ref{E1}). It means that separatrix cycles can be formed or destroyed only under a variation of the parameter $\alpha$ for fixed $k_{1}$ and $k_{2}.$

As is known \cite{4}, the character of stability of a sewed nondegenerate separatrix cycle is determined by the sign of its saddle quantities which are always positive in our case, when the saddles are inside the boundary of even regions $II,$ $IV,\ldots,$ $2K$ and $k_{2}>1$ (Theorems~44 and~47 from \cite{2} are valid for the piecewise linear dynamical systems as well). It follows that a nondegenerate separatrix cycle of (\ref{E1}) is always unstable (inside and outside) and, under a variation of $\alpha,$ such a separatrix cycle can generate at most $k+1$ small unstable limit cycles inside its loops (digons) or the only big unstable limit cycle outside it when system (\ref{E1})  has $2k+1$ simple singular points. The lemma is proved.$\qquad\square$

\section{Global Analysis}

For the global analysis of limit cycle bifurcations, we will use the Wintner--Perko termination principle
which is formulated for planar relatively prime analytic systems connecting the main local bifurcations of limit
cycles~\cite{8},~\cite{14}.

\medskip
\noindent\textbf{Theorem 5.1 (Wintner--Perko termination principle).}
    \emph{Any one-para\-me\-ter fa\-mi\-ly of multiplicity-$m$ limit cycles of a relatively prime analytic system 
\begin{equation}
	\label{E12}
    	\mbox{\boldmath$\dot{x}$}=\mbox{\boldmath$f$}
    	(\mbox{\boldmath$x$},\mbox{\boldmath$\mu$)},
\end{equation}
where $\mbox{\boldmath$x$}\in\textbf{R}^2;$ $\mbox{\boldmath$\mu$}\in\textbf{R}^n;$ $\mbox{\boldmath$f$}\in\textbf{R}^2$ $(\mbox{\boldmath$f$}$ is an analytic vector function), 
can be extended in a unique way to a maximal one-parameter family of multiplicity-$m$ limit cycles of
(\ref{E12}) which is either open or cyclic.}

\emph{If it is open, then it terminates either as the parameter or the limit cycles become unbounded; 
or, the family terminates either at a singular point of (\ref{E12}), which is typically a fine focus of 
multiplicity~$m,$ or on a $($compound) separatrix cycle of (\ref{E12}), which is also typically of
 multiplicity~$m.$}
\medskip

As was shown in \cite{4}, the return map constructed in a neighborhood of a multiple limit cycle of a piecewise linear dynamical system (under some general assumptions which are valid in our case) will be an analytic function; therefore, all statements formulated for multiplicity-$m$ limit cycles of analytic system (\ref{E12}) will be also valid for sewed multiplicity-$m$ limit cycles. Applying the Wintner--Perko termination principle and the results on local bifurcations of limit cycles obtained in the previous section for system (\ref{E1}), we will prove the following theorem. 

\medskip
\noindent \textbf{Theorem 5.2.}
\emph{System (\ref{E1}) with $k$ dropping sections and $2k+1$ singular points can have at most $k+2$ limit cycles, $k+1$ of which surround the foci one by one and the last, $(k+2)$-th, limit cycle surrounds all of the singular points of~(\ref{E1}).}
\medskip

\noindent \textbf{Proof.} Let us consider global bifurcations of limit cycles of (\ref{E1}).  Suppose again that the zero isocline $\beta-\alpha x-y=0$ passes through the corner point $(x_{1},y_{1})$ of the infinite isocline $y=\varphi(x)$ and that $\alpha>\alpha^{*}.$ In this case, the only singular point in the phase plane is a sewed stable focus and all trajectories of (\ref{E1}) tend to it when $t\rightarrow+\infty.$ For  decreasing $\alpha$ $(k_{2}<\alpha<\alpha^{*}),$ the sewed focus becomes unstable and a stable limit cycle is generated from the boundary curve of the domain with closed trajectories (immediately after passing the value $\alpha^{*}$ by the parameter $\alpha).$

For $\alpha=k_{2},$ the first dropping section of (\ref{E1})  will coincide with a part of the straight line $\beta-\alpha x-y=0$  and an unstable equilibrium segment will appear inside the stable limit cycle. If we rotate the line $\beta-\alpha x-y=0$ around an interior point of the segment (changing both of the parameters, $\alpha$ and $\beta),$ two unstable limit cycles surrounding stable foci one by one will be generated from the ends $(x_{1},y_{1})$ and $(x_{2},y_{2})$ of the equilibrium segment. Under the further rotation of the line $\beta-\alpha x-y=0,$ it will pass first through the next corner point, $(x_{4},y_{4}),$ and then, successively, through the points $(x_{6},y_{6}),\ldots,$ $(x_{2k},y_{2k}).$  Every time, the corner point becomes a sewed saddle-node generating an unstable limit cycle surrounding a stable focus. So, we will get a piecewise linear system with $2k+1$ singular points having at least $k+1$ small unstable limit cycles surrounding the stable foci one by one inside a big stable limit cycle, $k+2,$ surrounding all of the singular points.

Under the further rotation of the zero isocline, all $k+1$ small limit cycles simultaneously disappear in a nondegenerate separatrix cycle consisting of $k+1$ loops (digons), this separatrix cycle generates a big (unstable) limit cycle which combines with another big (stable) limit cycle of (\ref{E1}) forming a semi-stable (double) limit cycle which finally disappears in a so-called trajectory condensation.

Let us prove now that system (\ref{E1}) cannot have more than  $k+2$ limit cycles. The proof is carried out by contradiction. Since a small limit cycle is always unique in the corresponding strip regions (Lemma~4.1), suppose that system  (\ref{E1}) with three field rotation parameters, $k_{1},$ $k_{2},$ and $\alpha,$ has three big limit cycles. Then we get into some domain in the space of these parameters which is bounded by two fold bifurcation surfaces forming a cusp bifurcation surface of multiplicity-three limit cycles \cite{8},~\cite{14}.

The cor\-res\-pon\-ding maximal one-parameter family of multiplicity-three limit cycles cannot be cyclic, otherwise there
will be at least one point cor\-res\-pon\-ding to the limit cycle of multi\-pli\-ci\-ty four (or even higher) in the parameter
space. Extending the bifurcation curve of multi\-pli\-ci\-ty-four limit cycles through this point and parameterizing the
corresponding maximal one-parameter family of multi\-pli\-ci\-ty-four limit cycles by a field rotation para\-me\-ter,
for example, by the parameter~$\alpha,$ we will obtain a monotonic curve (see \cite{8}, \cite{14}) which, by the Wintner--Perko termination principle (Theorem~5.1), terminates either at the boundary curve of the domain filled by closed trajectories of (\ref{E1}) or on some degenerate separatrix cycle of (\ref{E1}).  Since we know at least the cyclicity of the boundary curve which is equal to two (Lemma~4.2), we have got a contradiction with the termination principle (Theorem~5.1) stating that the multiplicity of limit cycles cannot be higher than the multi\-pli\-ci\-ty (cyclicity) of the end bifurcation points in which they terminate.

If the maximal one-parameter family of multiplicity-three limit cycles is not cyclic, using the same principle (Theorem~5.1), this
again contradicts with the cyclicity result for the boundary curve (Lemma~4.2) not admitting the multiplicity of limit cycles higher than two. Moreover, it also follows from the termination principle that the degenerate separatrix cycles of (\ref{E1}) cannot have the multiplicity (cyclicity) to be higher than two. Therefore, according to the same principle, there are no more than two big limit cycles in the exterior domain outside the boundary curve of (\ref{E1}). 

The same results can be obtained by means of the new geometric methods developed in \cite{9},  \cite{10}. The phase portraits and bifurcation diagrams for system (\ref{E1}) will be similar to that which were constructed in \cite{3}, \cite{4}. Thus, system (\ref{E1}) with $2k+1$ singular points cannot have more than $k+2$ limit cycles, i.\,e., $k+2$ is the maximum number of limit cycles of such system and the obtained distribution $(k+1$ small limit cycles plus a big limit cycle) is the only possibility for their distribution. This completes the proof of the theorem.$\qquad\square$

\end{document}